\documentclass[graybox]{svmult}

\usepackage{type1cm}        
\usepackage{makeidx}         
\usepackage{graphicx}        
\usepackage{multicol}        
\usepackage[bottom]{footmisc}

\usepackage{newtxtext}       %
\usepackage{newtxmath}       

\usepackage[T1]{fontenc}
\usepackage[hidelinks]{hyperref}
\usepackage[ruled, vlined, linesnumbered]{algorithm2e}
\usepackage{stmaryrd}

\usepackage{codehighlight}
\usepackage{pifont}
\usepackage{listings}
\usepackage{caption}

\usepackage{amsopn}
\DeclareMathOperator{\diag}{diag}

\newcommand{\pluseq}{\mathrel{+}=}
\newcommand*{\Scale}[2][4]{\scalebox{#1}{$#2$}}
\newcommand{\dual}{^{\prime}}

\newcommand{\reals}{{\mathbb{R}}}

\newcommand{\semi}[1]{\lvert#1\rvert}

\makeindex             

\begin{document}

\title*{Assembly of multiscale linear PDE operators}

\author{Miroslav Kuchta}
\institute{Miroslav Kuchta \at Simula Research Laboratory, P.O. Box 134, 1325 Lysaker, Norway, \email{miroslav@simula.no}}
%
%
\maketitle


\abstract{
  In numerous applications the mathematical model consists of different
  processes coupled across a lower dimensional manifold. Due to the multiscale
  coupling, finite element discretization of such models presents a challenge.
  Assuming that only \emph{singlescale} finite element forms can
  be assembled we present here a simple algorithm for representing
  multiscale models as linear operators suitable for Krylov methods. Flexibility
  of the approach is demonstrated by numerical examples with coupling across
  dimensionality gap 1 and 2. Preconditioners for several of the problems are discussed.
}

\section{Introduction}\label{sec:intro}
This paper is concerned with implementation of the finite element method (FEM) for 
multiscale models, that is, systems where the unknowns are defined over domains of (in general)
different topological dimension and are coupled on a manifold, which is possibly
a different domain. 
The systems arise naturally in applications where Lagrange multipliers are used to
enforce boundary conditions, e.g. \cite{babuvska1973finite, bertoluzza2017boundary}, 
or interface coupling conditions e.g. \cite{bernardi1993domain, ambartsumyan2018lagrange, layton2002coupling}. 
 In modeling reservoir flows \cite{cerroni2019mathematical}, tissue perfusion
\cite{cattaneo2014computational, d2008coupling, koch2019modeling} or soil-root interaction \cite{koch2018new}
resolving the interface as a manifold
of co-dimension 1 can be prohibitively expensive. In this case it is convenient to represent
the three-dimensional structures as curves and the model reduction gives rise to multiscale
systems with a dimesionality gap 2.

Crucial for the FEM discretization of the multiscale models is the assembly of
coupling terms, in particular, integration over the coupling manifold. There exists
a number of open source FEM libraries, e.g. \cite{BangerthHartmannKanschat2007, mfem-library, MR3043640, getfem}, which expose this
(low-level) functionality and as such can be used for implementation.
However, for rapid prototyping, it is advantageous if the new models are described in a more abstract
way which is closer to the mathematical definition of the problem.

FEniCS is a popular open source FEM framework which employs a compiler 
to generate low level (C++) assembly code from the high-level symbolic representation of the variational
forms in the UFL language
embedded in Python, see \cite{logg2012automated}. Here the code generation pipeline provides
convenience for the user. At the same time, implementing new features is complicated by the fact
that interaction with all the components of the pipeline is required. 
As a result, support for multiscale models has only recently been added to the core
of the library \cite{daversin2019abstractions} and is currently
limited to problems with dimensionality gap 0 and 1. Moreover, in case of the
trace constrained systems the coupling manifold needs to be
triangulated in terms of facets of the bulk discretization. We remark that
similar functionality for multiscale systems is offered by the FEniCS based
library \cite{multifenics}.

Here we present a simple algorithm\footnote{
  Implementation can be found in the Python module \texttt{FEniCS\textsubscript{ii}} \url{https://github.com/MiroK/fenics_ii}}
which extends FEniCS to support a more
general class of multiscale systems by transforming symbolic variational
forms in UFL language into a domain specific language \cite{Mardal2012} which represents
(actions of) discrete linear operators. As this representation targets solutions
by iterative methods preconditioning strategies shall also be discussed. 
Our work is structured as follows. Section \ref{sec:algorithm} details the algorithm. 
Numerical examples spanning dimensionality gap 0, 1 and 2 are presented in 
\S\ref{sec:trace} and \S\ref{sec:multi} respectively.

\section{Multiscale assembler}\label{sec:algorithm}
In the following $(\cdot, \cdot)_{\Omega}$ denotes the $L^2$ inner
product over a bounded domain $\Omega\subset\reals^{d}$, $d=1, 2, 3$. The duality
pairing between the Hilbert space $V$ and its dual space $V\dual$ is denoted by
$(\cdot, \cdot)$. Given basis of a discrete finite element space $V_h$,
the matrix representation of operator $A$ is $A_h$. Adjoints of $A$ and $A_h$ 
are denoted as $A\dual$ and $A_h\dual$ respectively.

Our representation of multiscale systems builds on two observations, which
shall be presented using the Babu{\v s}ka problem \cite{babuvska1973finite}. Let 
$\Gamma=\partial\Omega$ and $V=H^1(\Omega)$, $Q=H^{-1/2}(\Gamma)$, $W=V\times Q$. 
Then for every $L\in W\dual$ there exists a unique solution $w=(u, p)\in W$ satisfying 
$\mathcal{A}w=L$ where
\begin{equation}\label{eq:babuska}
  \mathcal{A}=\begin{pmatrix}
    A & B\dual\\
    B & 0
  \end{pmatrix}
  \quad\mbox{and}\quad
  \begin{aligned}
     (Au, v) &= (\nabla u, \nabla v)_{\Omega} + (u, v)_{\Omega}\quad v\in V,\\
     (Bu, q) &= (Tu, q)_{\Gamma}\quad q \in Q.
  \end{aligned}
\end{equation} 
Here $T:H^1(\Omega)\rightarrow H^{1/2}(\Gamma)$ is the trace operator such that 
$Tu=u|_{\Gamma}$, $u\in C(\overline{\Omega})$. We remark that \eqref{eq:babuska} is the 
weak form of $-\Delta u + Iu = f$ in $\Omega$ with $u=g$ on $\partial\Omega$ enforced 
by the Lagrange multiplier $p$.

Given the structure of $\mathcal{A}$ in \eqref{eq:babuska} it is natural to represent 
the operator on a finite element space $W_h$ as a block structured matrix (rather then 
a monolithic one). Moreover, observe that the \emph{multiscale} operator $B: V\rightarrow Q$ 
(operator $A: V\rightarrow V\dual$ is \emph{singlescale}) is a composition $B=I\circ T$ where 
$I:H^{1/2}(\Gamma)\rightarrow Q$ is a singlescale operator. Therefore, matrix representation of 
$B$ is a matrix product $B_h=I_h T_h$. Assuming that the FEM library at hand can only assemble 
singlescale operators, e.g. $I$ and $A$, the multiscale operators $B_h$ and $\mathcal{A}_h$ can be formed if 
representation of the trace operator is available. We remark that the block representation
is advantageous for construction of preconditioners; for example the blocks can be easily
shared between the system and the preconditioner, cf. \cite{kirby2018solver, Mardal2012}.

Based on the above observations the multiscale systems can be represented as block structured 
operators where the blocks are not necessarily matrices. Cbc.block \cite{Mardal2012} defines a 
language for matrix expressions using the lazy evaluation pattern. In particular, 
block matrix(\texttt{block\textunderscore mat}) and matrix product(\texttt{\textasteriskcentered}) are built-in 
operators. We remark that the operators are not formed explicitly, however, they can be evaluated 
if e.g. action in a matrix-vector product in a Krylov solver is needed. Using $B$ from \eqref{eq:babuska}
as an example we thus aim to build an interpreter which translates UFL representation of $(Tu, q)_{\Gamma}$ 
into a cbc.block representation $I_h \text{\textasteriskcentered} T_h$. We remark that $T_h$ is here
assumed to be a mapping between primal representations, cf. \cite{mardal2011preconditioning}.

The core of the multiscale interpreter is the algorithm (Figure \ref{fig:algos})
translating between the two symbolic representations. Observe that in \texttt{multi\textunderscore assemble} different
\emph{reduced} assemblers are recursively called on the transformed UFL form with the singlescale
form being the base case. An example of a reduced assembler is the \texttt{trace\textunderscore assemble}
function which, having found \emph{trace}-reduced argument (ln. \ref{terminal}) in form $a$, e.g. $a(u, q)=(Bu, q)=(Tu, q)_{\Gamma}$,
$u\in V_h$, $q\in Q_h$ builds a finite element \emph{trace space} 
$\bar{V}_h=\bar{V}_h(\Gamma)$ (ln. \ref{trace-space}), an algebraic representation of the operator
$T:V_h\rightarrow \bar{V}_h$ (ln. \ref{trace-mat}) and delegates
assembly of the transformed form $I(\bar{u}, q)=(\bar{u}, q)_{\Gamma}$, $\bar{u}\in \bar{V}_h$,
$q\in Q_h$ (ln. \ref{reduced-form}) to \texttt{multi\textunderscore assemble} (ln. \ref{multi-pass}). As $I$ is singlescale the
native FEniCS assemble function can be used to form the matrix $I_h$ and
the symbolic matrix-matrix product representation can be formed (ln. \ref{multi-pass}).
The translation can thus be summarized as $(Tu, q)_{\Gamma}\rightarrow (\bar{u}, q)_{\Gamma}\text{\textasteriskcentered}T_h\rightarrow I_h \text{\textasteriskcentered} T_h$.

\begin{figure}[]
\begin{minipage}{0.465\textwidth}
    \begin{algorithm}[H]\scriptsize{
      \KwData{a::UFL.Form or list of UFL.Form}
      \KwResult{cbc.block matrix expression}
      \Begin{
        \tcp{Single form}
        \If{a is UFL.Form}{
          \tcp{Attempt to reduce}
          \For{assemble $\in$ assemblers}{
            tensor = assemble(a)\\
            \If{tensor is not \normalfont{\textbf{None}}}{
              \KwRet{tensor}
            }    
          }
          \tcp{Singlescale operator}
          \KwRet{FEniCS.assemble(a)}
        }
        \tcp{Functional}
        \If{is\textunderscore number(form)}{\KwRet{form}}
        \BlankLine
        shape $\leftarrow$ sizes(form)\\
        \tcp{Assemble blocks}
        blocks $\leftarrow$ \textbf{map}(multi\textunderscore assemble, form)\\
        \tcp{List/List of operators}
        tensor $\leftarrow$ reshape(blocks, shape)\\
        \BlankLine
        \tcp{Reshape for cbc.block}
        \tcp{Form had test functions only}
        \If{is\textunderscore vector(tensor)}{
          \KwRet{block.block\textunderscore vec}(tensor)
        }
        \tcp{Bilinear form}
        \KwRet{block.block\textunderscore mat}(tensor)
      }
      \caption{multi\textunderscore assemble\label{algo:multi}}
    }
    \end{algorithm}
\end{minipage}
\begin{minipage}{0.535\textwidth}
    \begin{algorithm}[H]\scriptsize{
      \KwData{a::UFL.Form}
      \KwResult{cbc.block matrix expression}
      \Begin{
        trace\textunderscore integrals $\leftarrow$ get\textunderscore trace\textunderscore integrals(a)\\
        all\textunderscore integrals $\leftarrow$ integrals(form)\\
        \If{\textnormal{\textbf{not}} trace\textunderscore integrals}{
          \KwRet{\textnormal{\textbf{None}}}
        }
        \tcp{Form is sum of integrals...}
        cs $\leftarrow$ []\\
        \For{i $\in$ all\textunderscore integrals}{
          \If{i $\notin$ trace\textunderscore integrals}{
            cs $\pluseq$ [multi\textunderscore assemble(Form([i]))]\\
            \textbf{continue}
          }
          intgrnd $\leftarrow$ integrand(i)\\ 
          $u, \leftarrow$ trace\textunderscore terminals(intgrnd)\nllabel{terminal}\\
          $V_h \leftarrow$ function\textunderscore space($u$)\\
          $\bar{V}_h \leftarrow$ trace\textunderscore space($V_h, u$)\nllabel{trace-space}\\
          $T_h \leftarrow$ trace\textunderscore matrix($V_h$, $\bar{V}_h$)\nllabel{trace-mat}\\
          
          \If{is\textunderscore trial\textunderscore function(u)}{
            $\bar{u} \leftarrow$ TrialFunction($\bar{V}_h$)\\
            ii $\leftarrow$ replace(intgrnd, u, $\bar{u}$)\\
            $I$ = Form([reconstruct(i, ii)])\nllabel{reduced-form}\\
            $B_h\leftarrow$ multi\textunderscore assemble($I$)\textasteriskcentered T\nllabel{multi-pass}\\
            cs $\pluseq$ [$B_h$]
          }

          \tcp{Handle test/function}
        }
        \tcp{...cbc.block sum of operators}
        \KwRet{\textnormal{\textbf{reduce}}(+, cs)}\nllabel{product}
      }
    \caption{trace\textunderscore assemble\label{algo:trace}}
    }
    \end{algorithm}
    \end{minipage}
  \vspace{-10pt}        
\caption{Translation of UFL representation of multiscale variational form into cbc.block 
matrix expression. Several passes by different scale assemblers might be needed to 
reduce the form into base case singlescale which can be assembled as matrix or vector by FEniCS.
Handling of test function and function type terminals is omitted for brevity.
}
\label{fig:algos}
\end{figure}

Algorithm \ref{algo:multi} can be easily extended to different multiscale couplings
by adding a dedicated assembler. In particular, given $\Omega\subset\reals^3$ and $\gamma$
a curve contained in $\Omega$, the 3$d$-1$d$ coupled problems \cite{d2008coupling, cerroni2019mathematical}
require operators $T$, $\Pi$ such that for $u=C(\Omega)$, $Tu=u|_{\gamma}$  and
\begin{equation}\label{eq:Pi_avg}
(\Pi u)(x) = \lvert C_R(x) \rvert ^{-1}\int_{C_R(x)} u(y)\,\mathrm{d}y.
\end{equation}
Here $C_R(x)$ is a circle of radius $R$ in a plane
$\{y\in\mathbb{R}^3, (y-x)\cdot\tfrac{\mathrm{d}\gamma}{\mathrm{d} s}(x) = 0\}$
defined by the tangent vector of $\gamma$ at $x$. We remark that assembling 3$d$-1$d$ constrained
operators follows closely Algorithm \ref{algo:trace}, with the non-trivial difference being
the representation of $\Pi$. We remark that in assembly of $\Pi$ or $T$ we do not
require that $\gamma$ is discretized in terms of edges of the mesh of $\Omega$. In fact,
the two meshes can be independent. This is also the case for $d$--$(d-1)$ trace.
Let us also note that the restriction operator $Ru=u|_{\omega}$, where $\omega\subseteq\Omega\subset\reals^d$
can be implemented similar to the trace operator. Finally, observe that the Algorithm
\ref{algo:multi} is not limited to forms where the arguments are \emph{reduced} to the
coupling manifold. Indeed, \cite{cerroni2019mathematical, gjerde2018singularity} utilize
\emph{extension} from $\gamma$ to $\Omega$ by a constant or as Green function of a line source respectively.
Such couplings can be readily handled if realization of the discrete \emph{extension} operator is
available.

We conclude the discussion by listing the limitations of our current implementation.
Unlike in \cite{daversin2019abstractions, multifenics} the MPI-parallelism is missing\footnote{
  The serial performance of our pure Python implementation is cca. 2x slower
  than the native FEniCS implementation \cite{daversin2019abstractions}. More precisely, assembling \eqref{eq:babuska} on
  $\Omega=\left[0, 1\right]^2$ discretized by $2\cdot 1024^{2}$ triangles and continuous linear
  Lagrange elements (the system matrix size is approx $10^6$, however, it is \emph{not} explicitly formed here)
  takes 3.86s (to be compared with 1.79s). Most of the time is spent building $T_h$.
  The trace matrix is reused by the interpreter to evaluate both $B_h$ and $B\dual_h$.
} as is the support for nonlinear forms. Moreover, the reduction operators cannot be
nested and can only be applied to terminal expressions in UFL, e.g. $T(u+v)$ cannot be interpreted.
In addition, point constraints are not supported. With the exception of parallelism the
limitations should be addressed by future versions.

In the following we showcase the multiscale interpreter by considering coupled
problems with dimensionality gap 0, 1 and 2. We begin by a trace constrained
2$d$-1$d$ Darcy-Stokes system. 

\section{Trace constrained systems}\label{sec:trace}
Let $\Omega_1$, $\Omega_2\subset\reals^2$ be such that $\Gamma=\partial{\Omega_1}\cap\partial{\Omega_2}$ and $\semi{\Gamma}\neq 0$.
Further let $\partial\Omega_i=\Gamma\cup\Gamma^{D}_i\cup\Gamma^{N}_i$ where $\semi{\Gamma^{k}_i}\neq 0$, $i=1, 2$, $k=N, D$
and $\Gamma\cap\Gamma^{N}_{i}=\emptyset$, cf.Figure \ref{fig:DS_domain}. We then
wish to solve the Darcy-Stokes problem (with unit parameters)

\begin{minipage}{0.35\textwidth}
  \centering
  \includegraphics[width=\textwidth]{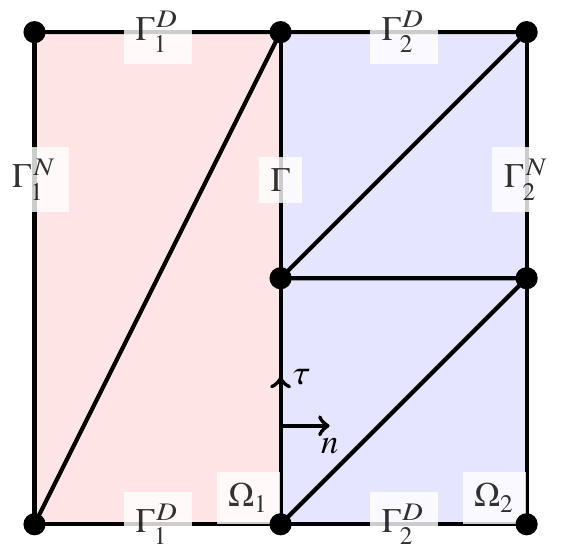}
  \vspace{-25pt}
    \captionof{figure}{Domain for \eqref{eq:darcy_stokes_strong}.}
    \label{fig:DS_domain}
\end{minipage}
\begin{minipage}{0.6\textwidth}
  \centering
  {
  \begin{subequations}\label{eq:darcy_stokes_strong}
    \begin{align}
      \label{eq:darcystokes1} 
      -\nabla\cdot\sigma &= f_1  &\text{ in } \Omega_1,\\
      \label{eq:darcystokes2}
      \nabla \cdot u_1 &= 0 &\text{ in } \Omega_1,  \\
      \label{eq:darcystokes3} 
      u_2 + \nabla p_2 &= 0  &\text{ in } \Omega_2,\\
      \label{eq:darcystokes4}
      \nabla \cdot u_2 &= f_2 &\text{ in } \Omega_2,  \\
      \label{eq:darcystokesmass} 
      u_1 \cdot n - u_2 \cdot n &= 0 &\text{ on } \Gamma, \\
      \label{eq:darcystokesstress}
      n \cdot\sigma\cdot n + p_2 &=0  &\text{ on } \Gamma ,\\
      \label{eq:darcystokesBJS}
  -n\cdot\sigma\cdot \tau - u_1 \cdot \tau &= 0  &\text{ on } \Gamma.
\end{align}
  \end{subequations}
  }
\end{minipage}
Here $\sigma(u_1, p_1)=D(u_1)-p_1I$ with $D(u)=\tfrac{1}{2}((\nabla u)+(\nabla u)\dual)$. The unknowns
$u_1$, $p_1$ and $u_2$, $p_2$ are respectively the Stokes and Darcy velocity and pressure.
The system is closed by prescribing Dirichlet conditions on $\Gamma^D_i$ and
Neumann conditions on $\Gamma^N_i$.

Let $T_n$, $T_t$ be the normal and tangential trace operators on $\Gamma$.
We shall consider variational formulations of \eqref{eq:darcy_stokes_strong} induced
by a pair of operators
\begin{equation}\label{eq:DS_op}
  \mathcal{A}_p=\left(\begin{array}{cc|c}
  -\nabla\cdot D + T\dual_tT_t & -\nabla & T_n\dual\\
  \text{div} & & \\
  \hline
  -T_n &       &                            -\Delta
  \end{array}\right),
  \mathcal{A}_m=\left(\begin{array}{cc|cc|c}
  -\nabla\cdot D + T\dual_tT_t & -\nabla  &   &           & T\dual_n\\
  \text{div}                  &         &    &          &          \\
  \hline
                               &         & I  & -\nabla & -T\dual_n\\
  &         &\text{div}  & &         \\
  \hline
  T_n                          &         & -T_n        & &\\
  \end{array}\right).
\end{equation}
Using the (mixed) operator $\mathcal{A}_m$ problem \eqref{eq:darcy_stokes_strong}
is solved for both $u_2$, $p_2$ and an additional unknown, the Lagrange multiplier,
which enforces mass conservation $u_1 \cdot n - u_2 \cdot n = 0$ on $\Gamma$. In
the (primal) operator $\mathcal{A}_p$ the condition appears naturally. Observe
that the operator is non-symmetric.

Well-posedness of the primal and mixed formulations as well the corresponding solution
strategies have been studied in a number of works, e.g \cite{discacciati2002mathematical}
and \cite{layton2002coupling, galvis2007non}. Here we compare the formulations
and discuss monolithic solvers which utilize block diagonal preconditioners
\begin{equation}\label{eq:DS_precond}
  \begin{aligned}
  \mathcal{B}_p=\diag\left(-\nabla\cdot D + T\dual_tT_t, I, -I\right)^{-1},\\
  \mathcal{B}_m=\diag\left(-\nabla\cdot D + T\dual_tT_t, I, I-\nabla\text{div}, I, (-\Delta+I)^{1/2}\right)^{-1}.
  \end{aligned}
\end{equation}
Here the preconditioner $\mathcal{B}_p$ has been proposed by \cite{cai2009preconditioning}, while 
$\mathcal{B}_m$ follows from the analysis \cite{galvis2007non} by operator preconditioning
technique \cite{mardal2011preconditioning}. More precisely, $\mathcal{B}_m$ is a Riesz map
with respect to the inner product of the space in which \cite{galvis2007non} prove well-posedness
of $\mathcal{A}_m$, i.e. $H^1_{\Scale[0.5]{{0, \Gamma^D_1}}}(\Omega_1)\times L^2(\Omega_1)\times H_{\Scale[0.5]{{0, \Gamma^D_2}}}(\text{div}, \Omega_2)\times L^2(\Omega_2)\times H^{1/2}(\Gamma)$.
We remark that all the blocks of the preconditioners can be realized by efficient and
order optimal multilevel methods. In particular, we shall use further the multigrid
realization of the fractional Laplace preconditioner \cite{baerland2018multigrid}.

In order to check mesh independence of the preconditioners let us consider
the geometry from Figure \ref{fig:DS_domain} and let $\Omega_1=\left[0, 0.5\right]\times\left[0, 1\right]$,
$\Omega_2=\left[0.5, 1\right]\times\left[0, 1\right]$. In both $\mathcal{A}_m$,
$\mathcal{A}_p$ the triangulations of the domains shall be \emph{independent}
\footnote{
  Details of experimental setup.
  We discretize $\Omega_i$ uniformly by first dividing the domains into $n\times m$ rectangles
  and afterwords splitting each rectangle into two triangles. For $\Omega_1$ we have $m=n$,
  $m=2n$ for $\Omega_2$ so that the trace meshes of the domains are different. Krylov solvers
  are started from random initial guess. Convergence tolerance for relative preconditioned
  residual norm of $10^{-10}$ is used. Unless specified otherwise the preconditioner blocks
  use LU factorization.
},
cf. Figure \ref{fig:DS_domain}, with the mesh of $\Gamma$ defined in terms of
facets of $\Omega_2$. Finally, the finite element approximation of $\mathcal{A}_p$
shall be constructed using $P_2$-$P_1$-$P_2$ elements\footnote{
Finite element space of continuous Lagrange
elements of order $k$ is denoted by $P_k$ while $RT_0$ denotes the space
of lowest order Raviart-Thomas elements.
  } while $P_2$-$P_1$-$RT_0$-$P_0$-$P_0$
is used for the mixed formulation $\mathcal{A}_m$.

  \begin{minipage}{\textwidth}
  \begin{minipage}[b]{0.50\textwidth}
    \centering
    \includegraphics[width=\textwidth]{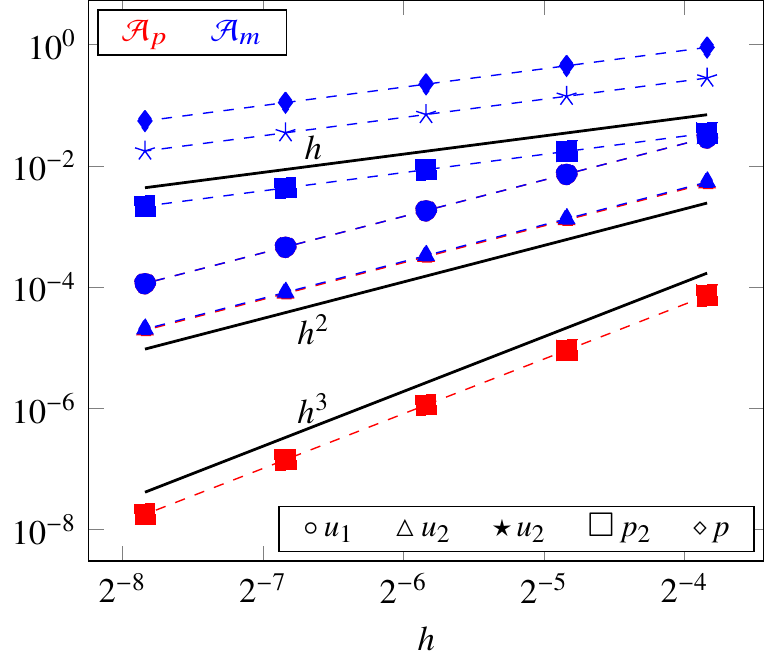}
    \label{fig:DS}
  \vspace{-25pt}    
  \captionof{figure}{
    Convergence of the primal(red) and mixed formulation of \eqref{eq:darcy_stokes_strong}.
    The approximation error is computed in the norms of $\mathcal{B}^{-1}_p$ and
    $\mathcal{B}^{-1}_m$.
  }
  \end{minipage}
  \vspace{5pt}
  \begin{minipage}[t]{0.45\textwidth}
    \centering
    \vspace{-190pt}
    \footnotesize{
      \captionof{table}{
        Number of iterations required for convergence of GMRes($\mathcal{A}_p$) and
        MinRes($\mathcal{A}_m$) using preconditioners \eqref{eq:DS_precond}, see also
        implementation in Figure \ref{code:DS}. Multigrid
        preconditioner for $H^{1/2}$ leads to slightly increased number
        of iterations compared to eigenvalue realization \cite{kuchta2016preconditioners}.
      }
      \label{tab:DS}      
      \begin{tabular}{c|ccc}
      \hline
      $h$ & $\mathcal{B}_p\mathcal{A}_p$ & $\mathcal{B}^{\Scale[0.5]{\text{EIG}}}_p\mathcal{A}_p$ & $\mathcal{B}^{\Scale[0.5]{\text{MG}}}_p\mathcal{A}_p$\\
      \hline
      $2^{-3}$ & 48 & 53 & 59 \\
      $2^{-4}$ & 48 & 51 & 59 \\
      $2^{-5}$ & 47 & 50 & 63 \\
      $2^{-6}$ & 47 & 49 & 65 \\
      $2^{-7}$ & 46 & 49 & 65 \\
      \hline
    \end{tabular}
    }
    \end{minipage}
  \end{minipage}

  Results of the numerical experiment are summarized in Table \ref{tab:DS}.
  It can be seen that the preconditioners \eqref{eq:DS_precond} are robust
  with respect to the discretization. Further, Figure \ref{fig:DS} shows
  that both formulations lead to expected order of convergence in all the unknowns.
  The approximation of Stokes variables is practically identical. We remark that
  $p_2$ convergence in $\mathcal{A}_p$ is reported in the $L^2$ norm for the sake
  of comparison with the mixed formulation. Implementation of $\mathcal{A}_m$ and
  preconditioner $\mathcal{B}_m$ can be found in Figure \ref{code:DS}.

  \begin{figure}[t]
    \begin{minipage}[H]{0.48\textwidth}
      \begin{python}
def mixed_darcy_stokes_system(n):
  '''Coupled Stokes-Darcy'''
  # Omega1 [0, 0.5]x[0, 1] as nxn, Omega2 nx2n
  # W = [P2]^2 x P1 x RT0 x P0 x P0
  W = [V1, Q1, V2, Q2, Q]

  u1, p1, u2, p2, p = map(TrialFunction, W)
  v1, q1, v2, q2, q = map(TestFunction, W)
  # Stokes traces
  Tu1, Tv1 = Trace(u1, gamma), Trace(v1, gamma)
  # Darcy traces
  Tu2, Tv2 = Trace(u2, gamma), Trace(v2, gamma)
  # Coupled integration
  dl = Measure('dx', domain=gamma)
  n, tau = Constant((1, 0)), Constant((0, 1))
  
  a = block_form(W, 2)
  # Stokes
  a.add(inner(sym(grad((u1)), sym(grad(v1)))*dx +
        inner(dot(Tu1, tau), dot(Tv1, tau))*dl
        -inner(q1, div(u1))*dx
        -inner(p1, div(v1))*dx)
  # Darcy
  a.add(inner(u2, v2)*dx-inner(p2, div(v2))*dx-
        inner(q2, div(u2))*dx)
  # Coupling
  a.add(
   inner(p, dot(Tv1, n))*dl-inner(p, dot(Tv2, n))*dl
   +inner(q, dot(Tu1, n))*dl-inner(q, dot(Tu2, n))*dl
  )
  # Define rhs + boundary conditions
  A, b = map(ii_assemble, (a, L))

  return A, b, W
  \end{python}
    \end{minipage}
    \hspace{5pt}
    \begin{minipage}[H]{0.45\textwidth}
      \begin{python}
def mixed_darcy_stokes_preconditioner(W, AA):
  '''H1 x L2 x Hdiv x L2 x H^{0.5}'''
  V1, Q1, V2, Q2, Q = W
  # Stokes velocity
  V1r = LU(AA[0][0])
  # Stokes pressure
  p, q = TrialFunction(Q1), TestFunction(Q1)
  Q1r = LU(assemble(inner(p, q)*dx))  # Or AMG

  # Darcy velocity
  mesh2 = V2.mesh()
  bcs = DirichletBC(V2,
                    Constant((0, 0)),
                    'near(x[1]*(1-x[1]), 0)')

  u, v = TrialFunction(V2), TestFunction(V2)
  a = inner(u, v)*dx + inner(div(u), div(v))*dx
  L = inner(Constant((0, 0)), v)*dx
  # Need symmetric assembly
  Hdiv_inner, _ = assemble_system(a, L, bcs)
  V2r = LU(Hdiv_inner)  # or HypreAMS

  # Darcy pressure
  p, q = TrialFunction(Q2), TestFunction(Q2)
  Q2r = LU(assemble(inner(p, q)*dx))  # or AMG
  # Multiplier H^s norm by Eigvp ...
  # Qr = HsNorm(Q, s=0.5, bcs=False)**-1
  # ... or multigrid
  Qr = HsNormMG(Q, s=0.5, bdry=None, s=0.5, 
                mg_params={'nlevels': 3,
                           'eta': 0.4,
                           'macro_size: 1'})

  return block_diag_mat([V1r, Q1r, V2r, Q2r, Qr])
      \end{python}
    \end{minipage}
  \vspace{-15pt}    
    \caption{Implementation of mixed Darcy-Stokes problem. (Left)
      Definition of the problem operator. (Right) Complete implementation
      of $\mathcal{B}_m$ preconditioner using either eigenvalue \cite{kuchta2016preconditioners} or
      multigrid \cite{baerland2018multigrid} realization of the fractional Laplacian.
    }
    \label{code:DS}
  \end{figure}

\section{More general multiscale systems}\label{sec:multi}
To show flexibility of the interpreter we finally consider a simple prototypical
3$d$-1$d$ coupled problem and an extended Darcy-Stokes problem with 2$d$-2$d$-1$d$
coupling. We will present both problems before discussing the results.

Let $\Omega\subset\reals^3$ be a bounded domain and let $\gamma$ be a curve embedded
in $\Omega$. Assuming $\gamma$ is a representation of the vasculature (e.g. as center lines)
parameterized by arc length coordinate $s$ a model of tissue \emph{perfusion} by
\cite{d2008coupling} is given as
\begin{equation}\label{eq:daq}
  \begin{aligned}
  -\nabla\cdot(k\nabla u) + \beta(\Pi u - p)\delta_{\gamma} &= 0\quad &\mbox{ in }\Omega,\\
  -\frac{\mathrm{d}}{\mathrm{d}s}(\hat{k}\frac{\mathrm{d}}{\mathrm{d}s}p) - \beta(\Pi u - p) &= 0\quad &\mbox{ on }\gamma.
  \end{aligned}
\end{equation}
Here $k$, $\hat{k}$ are the conductivities of the tissue and the vasculature, while
$\beta$ is the permeability. Observe that the exchange term is localized in $\Omega$
by the Dirac function $\delta_{\Gamma}$.

Let next $\Omega_i\subset\reals^d$, $d=2, 3$, $i=1, 2$ be the fluid domain and a porous
domain which share a common interface $\Gamma$. A model for \emph{transport} of a scalar
$\phi$ in such a medium $\Omega=\Omega_1\cup\Omega_2$ was recently analyzed by
\cite{baier}. Here we shall consider a simplified, linearized version of the system
\begin{equation}\label{eq:baier}
    \begin{aligned}
      -\nabla\cdot\sigma + g\phi &= f_1  &\text{ in } \Omega_1,\\
      \nabla \cdot u_1 &= 0 &\text{ in } \Omega_1,  \\
      u_2 + \nabla p_2 +g\phi &= 0  &\text{ in } \Omega_2,\\
      \nabla \cdot u_2 &= f_2 &\text{ in } \Omega_2,  \\
      -\Delta\phi + \nabla\cdot{fu_1} + \nabla\cdot{f u_2} &= 0 &\text{ in }\Omega,\\
\end{aligned}
\end{equation}
where $g$ and $f$ are given vector and scalar fields on $\Omega$. We remark that
\eqref{eq:baier} is considered with the interface conditions \eqref{eq:darcystokesmass}-\eqref{eq:darcystokesBJS}.

Compared to Babu{\v s}ka problem \eqref{eq:babuska} or Darcy-Stokes problem \eqref{eq:DS_op}
systems \eqref{eq:baier} and \eqref{eq:daq} introduce new multiscale couplings
\begin{equation}\label{eq:baier_op}
  \footnotesize{
  \mathcal{A}_p = 
  \begin{pmatrix}
    -k\Delta + T\dual \Pi & \beta T\dual\\
    -\beta\Pi                 & -\hat{k}\Delta + \beta I 
  \end{pmatrix},\,    
  \mathcal{A}_t=\left(\begin{array}{cc|cc|c|c}
  -\nabla\cdot D + T\dual_tT_t & -\nabla  &   &           & T\dual_n & R\dual_1\\
  \text{div}                  &         &    &          &          &\\
  \hline
                               &         & I  & -\nabla & -T\dual_n & R\dual_2\\
  &         &\text{div}  & &         & \\
  \hline
  T_n                          &         & -T_n        & & &\\
  \hline
  \text{div}\circ R_1 & & \text{div}\circ R_2 & & & -\Delta
  \end{array}\right).
  }
\end{equation}
Indeed, in the perfusion operator $\mathcal{A}_p$ the test functions in the bulk are
reduced to $\gamma$ by a 3$d$-1$d$ trace operator while $\Pi$ in \eqref{eq:Pi_avg}
is used for the trial functions. The transport operator $\mathcal{A}_t$ then
uses restriction operators $R_i\phi=\phi|_{\Omega_i}$, $i=1, 2$ for $\phi\in C(\Omega)$.
We remark that differently weighted Sobolev spaces are required in order for the 3$d$-1$d$
reduction operators to be well defined, see \cite{d2008coupling}. In particular, the trace
operator requires higher than $H^1$ regularity.

We test the abilities of the assembler by considering FEM discretization
of \eqref{eq:daq} in terms of $P_1$-$P_1$ elements while \eqref{eq:baier}
shall be discretized by $P_2$-$P_1$-$RT_0$-$P_0$-$P_0$-$P_2$. Here the setup for
\eqref{eq:daq} mirrors \S\ref{sec:trace}. However, to simplify the restriction
the meshes for $\Omega_1$ and $\Omega_2$ are not independent. Instead, they are
defined using the triangulation of $\Omega$. The perfusion problem is then
setup on a uniform discretization of $\left[0, 1\right]^3$ with $\gamma$
a straight line which, in general, is not aligned with the edges
of the mesh of $\Omega$.

\begin{figure}
  \includegraphics[width=0.5\textwidth]{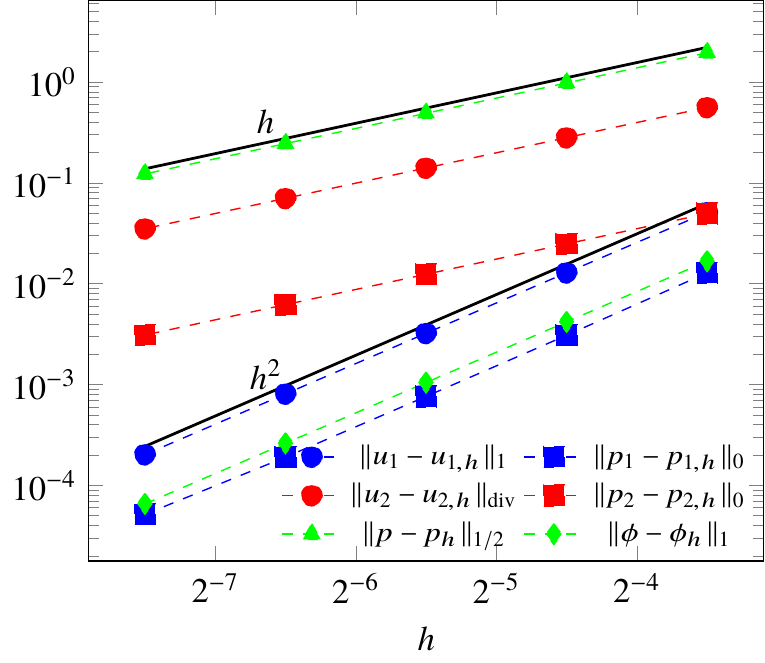}
  \includegraphics[width=0.5\textwidth]{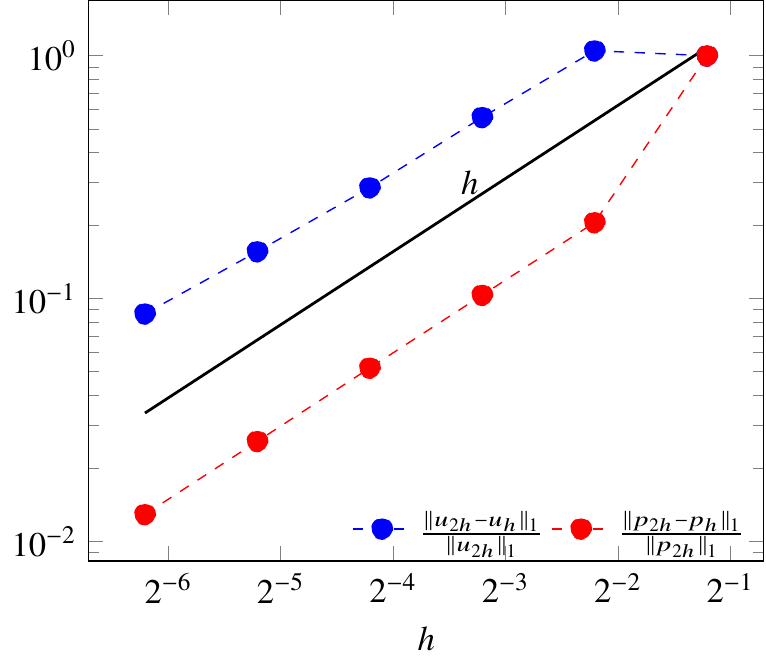}
  \vspace{-20pt}
  \caption{
    Convergence of the FEM approximation of the 2$d$-2$d$-1$d$ coupled problem
    \eqref{eq:baier} and a 3$d$-1$d$ problem \eqref{eq:daq}.
  }
  \label{fig:multi}
\end{figure}

Figure \ref{fig:multi} shows the error convergence of the two approximations.
For \eqref{eq:baier} the error with respect to the manufactured solution is
measured and the expected rates can be observed. In perfusion problem the relative
norm of the refined solution decreases linearly.

\bibliographystyle{spmpsci}
\bibliography{fenics_ii.bib}

\end{document}